\documentclass{article}
\pdfoutput=1

\usepackage{amssymb}

\usepackage{mathtools}
\newcommand{\no}{\nonumber}
\newcommand{\halb}{\frac{1}{2}}

\newcommand{\sI}{\mathbf{m}}
\newcommand{\Ncm}{N} 
\newcommand{\R}{\mathbb{R}} 
\newcommand{\eqcolon}{\mathrel{\resizebox{\widthof{$\mathord{=}$}}{\height}{ $\!\!=\!\!\resizebox{1.2\width}{0.8\height}{\raisebox{0.23ex}{$\mathop{:}$}}\!\!$ }}}
\DeclareMathOperator{\diag}{diag}

\title{Optimization-based modal decomposition for systems with multiple transports}

\author{Julius Reiss\thanks{Institut f\"ur Str\"omungsmechanik und Technische Akustik, Technical University Berlin, Germany, \texttt{reiss@tnt.tu-berlin.de}.}
}

\begin{document}

\maketitle

\begin{abstract}
Mode-based model-reduction is used to reduce the degrees of freedom of high dimensional systems, often by describing the system state by a linear combination of spatial modes.  
Transport dominated phenomena, ubiquitous in technical and scientific applications, often require a large number of linear modes to obtain a small representation error. 
This difficulty, even for the most simple transports,  originates from the inappropriateness of the decomposition structure in time dependent amplitudes of purely spatial modes.
     
In this article  an approach is discussed, which decomposes a flow field into  several fields of co-moving frames, where each one can be approximated by a few modes.
The method of decomposition  is formulated as an optimization problem. 
Different singular-value-based objective functions are discussed and connected to former formulations. 
A boundary treatment is provided. 
The decomposition is applied to generic cases and to a technically  relevant flow configuration of combustion physics.  
\end{abstract}

\noindent\emph{
keywords: 	
  model reduction, flow mechanics, non-linear model reduction, singular value decomposition, modal decomposition, boundary treatment, transport, POD \\
}
\noindent
AMS:	 65Z05, 35Q35

\section{Introduction}

Model reduction aims at reducing the computational complexity or the numerical cost of a systems by reducing its degrees of freedom. 
Often a mode-based approach is  used, since it allows to characterize the internal dynamics and  works for non-linear systems \cite{ChaS10}.   

In standard mode-based approaches, a snapshot-matrix, as a discrete representation of a field $q_{ij}= q(x_i,t_j)$, $q \in \R^{m,n}$ in space and time, is  approximated by
\begin{eqnarray}
q \approx  \tilde q  =   \sum_{l=1}^{r} \phi_l \otimes \alpha_l , \label{basicMode}  
\end{eqnarray}
 a sum of $r$  dyadic pairs $ (\phi_l \otimes \alpha_l)_{ij} =  \phi_l(x_i) \alpha_l(t_j)$ where  $\phi_l(x_i)$  are spatial modes and   $\alpha_l(t_j)$ are time dependent amplitudes.   
Parameter dependence is treated in the same way as  time dependence. 
The corresponding $r$-dimensional subspace can be used to build a reduced-order model (ROM), for instance by an interpolation procedure 
or by a Galerkin projection of the original system.

 Various decomposition methods  are known, such as the  balanced proper orthogonal decomposition \cite{TuR12,WilP02} or the dynamic mode decomposition 
\cite{CheTR12,SchmidSesterhenn2008,Tu13}. All share the structure of (\ref{basicMode}). 
One of the most popular modal decomposition approaches is the POD (proper orthogonal decompostion), \cite{HinV05,RatP03}, which technically reduces to a singular value decomposition (SVD) of the snapshot matrix.
The singular values depict the 2-norm associated with each  dyadic pair.
Thus,  a rapid singular value decay implies that a corresponding large part of the matrix is represented by the first pairs and, thereby,  a good low-rank approximation. 
The optimimality of the SVD in the 2-norm  is often referred to as   the Eckart-Young   theorem \cite{Stewart1991}. %
\smallskip

Transport-dominated phenomena usually pose a challenge for classical mode-based methods, since their dynamical behavior cannot be captured accurately by the linear combination of a few modes of the form \ref{basicMode}. 
A generic  example for  multiple  transports or waves is 
\begin{eqnarray}
q^{\mathrm{w}}_{ij}= q^w(x_i,t_j)=   \sum_{k=1}^{N} q_k(x_i -c_{k} t_j ,t_j ), 
\label{qw} 
\end{eqnarray}
where $q_k(x,t)$  are strongly localized functions in space $x$, and $c_k$ are wave velocities describing the translation with time $t_j$.  
We assume that functions $q_k(x,t)$  change slowly in time, i.e., is weakly dependent or even independent of its second argument. 
Physical realizations are  sound pulses, shock waves or flame fronts associated with rapid variation of pressure or temperature in space, and which are transported over long distances with a (nearly)  constant shape. 
The data $ q^{\mathrm{w}}$  typically needs large $r$ for a small approximation error $\| q-\tilde q\|$ in (\ref{basicMode}), in spite of its simple construction, 
which is  formalized  by  the  Kolmogorov $n$-width \cite{GreifUrban2019,OhlbergerRave2015}. 
As  example,  consider one wave ($N=1$)  with  $q_1$ and  $c_1$ chosen such that $q^w_{ik}=\delta_{i,j} $ is the (diagonal) unity matrix.  
Here, all singular values are one, so that no  decay is found which means that  a good approximation needs many modes \cite{RimMoeLeVeque2018}. 
As a problem of fundamental and practical importance it has recently attracted the attention of different groups suggesting various approaches to overcome this problem.

For  a single discontinuity the solution domain can be split by the position of the front
and a condition for the jump can be used \cite{TaddeiPerottoQuarteroni2015}, which resembles  classical shock fitting \cite{Moretti1987} of the numerical treatment of shocks.
A similar splitting of the domain is performed in \cite{ConstantineIaccarino2012},  where  the  Euler equations   are applied as the full model at the shock location. 
A dynamical procedure to handle moving shocks by splitting the support of the base functions to enrich the 
approximation subspace is presented in \cite{Carlberg2015}. 

Many of the transport dedicated  methods compensate the transports by a transformation, e.g. 
$T: r(x-ct,t ) \to r(x,t)  $,  and thereby  align the localized structure for all snapshots.   
This can strongly improve the convergence of the approximation \cite{Welper2017}.  
A time or parameter dependent shift can handle a single transport  
\cite{CagniartMadayStamm2019,LemkeMiedlarReissMehrmannSesterhenn2015,SesterhennShahirpour2016}. 
This is similar to the framework of symmetry reduction  \cite{RowleyKevrekidisMarsdenLust2003,FedeleAbessiRoberts2015}. 
Instead of a shift, a grid transformation can be used with a similar effect, as shown in one and two dimensions with complex shock and flame topologies  \cite{NairBalajewicz2019} or transport maps for parameter dependent  domains \cite{NoninoBallarinRozzaMaday2019}. 
The reduction of the case of a  flame front in two dimensions  with topology changes is reported in \cite{KrahSrokaReiss2019}.
To extend methods developed for a one dimensional to multidimensional setting, 
the Radon transformation can be used, producing a  set of one-dimensional projections \cite{Rim2018}. 
Moreover, the transports  challenge not only the reduction, but also the projections to create reduced models, which is discussed in 
\cite{AbgrallCrisovan2018,BrunkenSmetanaUrban2019,YanoPateraUrban2014}. 


To identify the  transports or transport velocities,  different approaches are available. 
If the transport is  considered  as the action of a Lie group, the correct transport 
is associated with minimal change of the transported quantity \cite{BeynThuemmler2004}, which was exploited for model reduction by \cite{OhlbergerRave2013}. 
 Transport velocities are identified by maximizing the  leading singular value as a function of the transport velocity in \cite{ReissSchulzeSesterhennMehrmann2018}. 
Optimal transport based on the concept of  Wasserstein distance is applied in \cite{IolloLombardi2014,BernardIolloRiffaud2018}, or similar 
 displacement interpolation \cite{RimMandli2018b} introduced by matching pieces of positive or negative slope. 
 An automatic velocity termination for the below discussed shifted POD was reported recently \cite{MendibleBruntonAravkinLowrieKutz2019}. 
 A formally different approach is  basis propagation by Lax-pair to create  time dependent bases \cite{GerbeauLombardi2014}. 
 
 A more automatic creation of reduced order models is obtained by the use of neural networks, where the transports are not explicitly handled, but 
 one relies on the general approximation strength of the nonlinear approach. 
 This can deliver very good reductions for the price of a reduced structural insight   \cite{LeeCarlberg2018,LuiWolf2019}.
\smallskip

In this report multiple transports, $N>1$ in \ref{qw}, are  considered. 
The central point is how to find, from given data $q^w$, a combination of $q_k$ which can be approximated by few modes.       
This is also addressed in  \cite{ReissSchulzeSesterhennMehrmann2018}, by  the author and co-workers,  as well as  by \cite{RimMoeLeVeque2018}.
Both methods use  transformations for multiple transports to align the transported quantities and reduce these with a truncated SVD.
The first method, the shifted POD, uses an iterative procedure to attribute components of the data to the velocity frame, where it gives lowest rank and finds,  for the tested simple hyperbolic examples ($q_k(x,t)\equiv r(x)$), the  expected representation by one dyadic pair  per transport  frame.  
In contrast, the second method uses a simple greedy approach, where a component associated with one frame cannot be re-attributed to a different, hindering a perfect decomposition. 
While the shifted POD  works well for the considered cases \cite{ReissSchulzeSesterhennMehrmann2018}, the given algorithm is difficult to analyze 
and the defining properties of determined decomposition are unclear. %

In this work, we (I) reformulate the decomposition, i.e., finding suitable $q_k$,  as an optimization problem based on the singular values in the different co-moving frames.
This not only allows to use various  methods of optimization,  but (II) provides  a better mathematical access to analyze and characterize the obtained results.  
We also seek to (III) reduce  non-uniqueness, which was observed in \cite{ReissSchulzeSesterhennMehrmann2018}.
Further, we (IV) introduce a simple and flexible boundary treatment.  
As in the previous method \cite{ReissSchulzeSesterhennMehrmann2018},  we (V) demand the  method to be purely data-based, since this promises a wide application area. 
We especially refrain from using special properties like hyperbolicity of the governing equation. 
First results of this approach were presented in \cite{Reiss2018}. 
We focus on the decomposition, while the construction of a dynamical reduced-order model on this base  is not within the scope of this report.
It is discussed for scalar laws e.g. in \cite{RimPeherstorferMandli2019} and for multiple transports in  \cite{BlackSchulzeUnger2019}, see also 
\cite{LeeCarlberg2018}. 

In \cite{SchulzeReissMehrmann2019}, the shifted POD decomposition was already formulated as an optimization, which is, in contrast  to the here proposed, based on the  residuum  of the approximation and not on the singular values. %
The methods are compared  in sec.~\ref{secFormer}.
\smallskip

The paper is organized as follows: 
In section \ref{secOpti} the optimization is formulated and the gradient of the objective functions is derived.
Uniqueness of the decompositions is discussed in sec. \ref{secExtensions}. 
The  method is  tested on  generic examples in sec.\ \ref{secFirstEx}. 
A generic  treatment of boundaries is given in sec.\ \ref{secBound}.
Sec.\ \ref{secInvestDecomp} discusses the resulting decomposition method and in sec.\ \ref{secFormer} the connection with the former  formulations of the shifted POD is given. 
In sec.\ \ref{secPDC}, the method is applied  to a reactive flow case,  showing a deflagration to detonation transition.  
Finally, we conclude and summarize in sec.\ \ref{secConcl}.

\section{Decomposition as an Optimization Problem}
\label{secOpti}

We aim for a  decomposition of a field  $q(x,t)$  into multiple transported quantities $q^k(x,t)$, 
assuming the field $q$ can be approximated by a structure like  (\ref{qw}).    
For transport dominated phenomena,   $q^k$ are expected to have  a more rapid singular value decay compared with $q$. 
A discrete setting   $(q)_{i,j}\coloneqq  q(x_i,t_j)$, $q\in \R^{m,n}$ is always assumed,  as typically  only discrete values from simulation or measurement are  available.  
Thus, we seek an representation  of the form  
 \begin{eqnarray}
  q &=&  \sum_{k=1}^{\Ncm}  T^{\Delta^k}\!\left[ q^{k}    \right] \label{decompGen}  \label{JRestrict}, 
 \end{eqnarray}
 with  the transformation operator
 as  a discrete approximation of the translation operation 
 $ 
 T^{\Delta^k}:   f(x,t)\to f(x+\Delta^k(t),t   )
 $.  
 $ T^{\Delta^k}[\cdot]$ maps data from the $\Ncm$ co-moving frames to the reference frame of the original data, called lab frame.
The form (\ref{decompGen})  accommodates problems like (\ref{qw}) with a time dependent shift  $\Delta^k_j= \Delta^k(t_j)$  presenting the location  of the  moving structures, $\Delta^k(t_j) = c_k t_j $ for the special case of constant transport velocities.
In this report we make  use of  the linearity  of the transformation $ T^{\Delta^k}[\cdot]$ in its argument $q^k$, while it is nonlinear in $\Delta^k_j$, since $q^k$ approximates a  nonlinear function $q^k(\cdot,t)$,  in general.  
Further, the existence of  an  inverse  $  T^{-\Delta^k} T^{\Delta^k} = T^{\Delta^k} T^{-\Delta^k} = \mathbf 1$ and orthogonality of  $  T^{\Delta^k}$  is assumed.    
These assumptions allow to replace a constrained by  an unconstrained  optimization,
 see  sec.~\ref{sec_grad}.

 The translation is  implemented   
by a Lagrange interpolation  (based on a fixed number of  neighboring points)  at the location shifted by $ \Delta^k_j$. 
Since shifting builds on  an interpolation, the transformation with negative shift might only approximate the inverse operation, as values at neighboring points are included in  $T^{-\Delta^k}[T^{\Delta^k}[\cdot]] $.

The relevance of (\ref{decompGen})  is that it allows to find a good approximation with smaller number of modes compared with  (\ref{basicMode})
 for  phenomena with multiple transports. %
Wave like phenomena  typically change slowly in the co-moving frame, so that $q^k$ has in this case  a much more rapid singular value decay than $q$. 
Thereby,   the residuum 
\begin{eqnarray}
R=   q-\tilde  q \label{residuum}
\end{eqnarray}
is expected to be  smaller
for the approximation  
$\tilde q =  \sum_{k=1}^{\Ncm}  T^{\Delta^k}\!\left[ \tilde q^{k}    \right] $ 
with 
 $\tilde q^k= \sum_{l=1}^{r_k}   \phi_{l}^k \otimes \alpha_{l}^k $, 
  compared with the form (\ref{basicMode}) with the same total number of modes.    
If the transports $\Delta^k$ are known or can be deduced from data, the  decomposition of the matrix $q$  into matrices  $ q^k$ remains as the  crucial step 
for a good reduction. 

\subsection{The Objective Function}
Finding the decomposition  (\ref{decompGen}), i.e., finding appropriate $q^k$  can be approached  in different ways.
One can try to identify structures which can be represented well in the different co-moving frames by reducing transformed data, as done in 
\cite{ReissSchulzeSesterhennMehrmann2018}.     
Or one   one can directly optimize  the norm of the residuum (\ref{residuum}), this approach was taken in \cite{SchulzeReissMehrmann2019}. 

Alternatively  one can %
define an arbitrary decomposition of $q \in \R^{m,n}$ into multiple co-moving frames with  values $q^k \in \R^{m,n}$ by using (\ref{JRestrict})  as a constraint.  
This ensures  that  the combination  of  the  data of the co-moving frames $q^k$  yields the original data.    
For any given decomposition, one can try to measure by an  appropriate expression the low rank approximation quality by the singular value decay of $q^k$. %
{Optimizing such an expression, or objective function,  by varying  the data  $q^k$ identifies the desired decomposition.  }


To  introduce a possible objective function for a given set of $q^k$,  the singular value decomposition    $  q^k = U^k S^k (V^k)^T $ is used, 
with the orthogonal matrices  $U^k \in  \R^{m,d}  $ and  $V^k\in \R^{n,d}$ and the diagnoal matrix $S^k \in \R^{d,d}$ containing the singular values 
 $ \diag(S^k) = s_1^k,s_2^k, \dots, s_d^k $ in descending order. 
Here,  $d = \min(m,n)$ is the maximal rank of the matrix. 
We define the  objective function  as  
\begin{eqnarray}
 J_2 &=& \sum_{k=1}^{\Ncm }  
\underbrace{
 \Bigl(1  - 
	 {
		\Bigl( 
		\sum_{l=1}^{r_k} \left(s^{k}_l\right)^2
		\Bigr)
      }
	/
	 { n^2_k}  
 \Bigr)
 }_{= J_2^k} 
  ,
\label{rankKJ}
\label{J2}
\end{eqnarray}    
with  the Frobenius norm $  n_k = \| q^k \|_F    $. 
Since  the Frobenius norm $n_k = \| q^k \|_F = \sqrt{\sum_{i,j} (q^k_{ij})^2 }$ satisfies  $n_k^2 = \sum_l (s^k_l)^2 $ and  $s^k_l\geq 0 $, it follows that $  n_k^2 \geq  \sum_{l=1}^{r_k} \left(s^{k}_l\right)^2  $,  so that $ J_2\geq 0 $ for any values $r_k\leq d$.    
If  the matrices $q^k$ are of rank $r_k$,  the sum of squares leading singular value equals the square of the norm and each term $J_2^k $ vanishes.
In this case,  $J_2$ attains its minimum zero. 
Minimizing (\ref{rankKJ}) respecting the constraint  (\ref{JRestrict}), equals optimizing towards  a 
decomposition into multiple matrices  of rank $r_k$.    
\medskip

This objective function is not the only possible way to measure the singular value decay. 
An alternative is to  build on the sum of the singular values, which  
is in analogy to the $l_1$ vector norm, which is used  to obtain sparse approximations:    
\begin{eqnarray}
J_1 &=& \sum_{k=1}^{\Ncm }  \left(  \sum_{l=1}^{d} s^k_l  \right)  = \sum_{k=1}^{\Ncm } \| q^k\|_*\label{J1}
\end{eqnarray}
This  1-norm of the singular values is called the Schatten 1-Norm $\| \cdot \|_* $ or nuclear norm or Ky-Fan norm. 
It was already noted by \cite{FazelHindiBoyd2001} that this serves as a heuristic rank minimization.

The objective function $J_1$ can be rewritten as the  Schatten-1-norm  of a single matrix by 
defining  a block matrix $Q$ composed from the matrices of the co-moving frames:
\begin{eqnarray}
 J_1 = \| Q\|_*  =  \left\| \begin{matrix}
	q^1 &  &&\\ 
	& q^2 &&\\ 
	&     & \ddots & \\
	&     &  & q^\Ncm \\
\end{matrix}  \right\|_*   \label{J1Schatten} 
\end{eqnarray} 
The last frame can be defined as 
$ 
q^\Ncm  =  T^{-\Delta^N}\left( q -\sum_{k=1}^{\Ncm -1}   T^{\Delta^k} [ q^k] \right) 
$
to include the constraint (\ref{JRestrict}), while all other  $q^k$ are arbitrary.

By rewriting $J_1$ as the   norm of a single matrix, the decomposition problem is reduced to the known problem of Schatten-1 optimization  
\cite{FazelHindiBoyd2001}.   %
Thereby, it is known that  $\min_{Q} J_1 $ is convex if the  set of all $Q$ is convex. 
Since the matrices $q^k,\, k=1,\dots, \Ncm-1$  are arbitrary and  $q^\Ncm$ is created from an affine map this is fulfilled.  
However, it is not strictly convex so that the minimum is not unique, but every minimum is a global minimum.\footnote{
	Remark: For Schatten-1-norm optimizations, dedicated methods are available, see e.g.\ \cite{LiuVandenberghe2010}.
Application of those  methods might offer methodological and practical benefits. 
Investigation of these methods is beyond the  scope of this report, concerned with \emph{finding} possible {formulations} of the decomposition problem. 
}

The disadvantage of $J_1$ is that it inherently needs  all singular values, which is    unfeasible for (typical) large problems. 
To mitigate this, we estimate the remaining singular values if only the first $r_k$ singular values are known. 
The total square sum  of the remaining singular values in each frame can be calculated from the Frobenius norm 
\begin{eqnarray}  
\sum_{l=r_k+1}^d  (s_l^k)^2 = n^2_k - \sum_{l=1}^{r_k}  (s_l^k)^2 = n^2_k J_2^k
\eqcolon \bar J_2^k.
\end{eqnarray}
The upper estimate for $J_1$ is obtained  for an equal distribution of $\bar J_2^k$ over all remaining singular values so that we approximate $\tilde s_l^k= \sqrt{\bar J_2^k/(d-r) } $ for $ l>r$, yielding   
\begin{eqnarray}
J_{1,2} &=& \sum_{k=1}^{\Ncm }  \left(  \sum_{l=1}^{r_k} s^k_l  + \sqrt{   (d-r)}\sqrt{  \bar J^k_2   }   \right) \label{J1delta}
\end{eqnarray}  
The second term   
$ \bar J_2^k %
$
equals up to a factor $n^2_k$ the components $J_2^k$ of the  objective function (\ref{rankKJ}). 
This suggests an alternative objective function 
 \begin{eqnarray}
\bar J_2 &=& \sum_{k=1}^{\Ncm }  
 	\Bigl(	{ n^2_k}  - 
 	{
 		\sum_{l=1}^{r_k} \left(s^{k}_l\right)^2
 	}  
 	\Bigr)
 	= 
 	\sum_{k=1}^{\Ncm } \bar J_2^k
 .
 \label{barJ2}
 \end{eqnarray}     
The difference between $J_2$ and $\bar J_2$ is discussed  further in section \ref{secInvestDecomp}.

\subsection{Deriving the Gradient} 
\label{sec_grad} 
For solving the optimization problem,  gradient based methods are used in this work. 
The gradient is  derived in two steps. 
First, the change of the objective function is evaluated per frame, second, the constraint (\ref{JRestrict})  is incorporated.

To this end, consider an (infinitesimal) small perturbation  of  the matrix entries  $ q^k_{ij} \to q^k_{ij}  + \delta q^k_{ij} $.  
For the Frobenius norm $  n^2_k = \sum_{i=1}^{m}\sum_{j=1}^{n} (q^k_{ij})^2$  the variation directly yields    
\begin{eqnarray}
 \frac{\partial (n^2_k) }{\partial  q^k_{ij}}  \,\delta q^k_{ij}  =  2  q^k_{ij}  \,\delta q^k_{ij}\;.
\end{eqnarray}
Assuming the singular values  to be simple,  the change of the $l^{th}$ singular value of $A= USV^T$  by $ A\to A+\delta A$ 
is 
$
u^T_l \delta A\, v_l = \delta s_l, \cite{Stewart1991}. 
$
Choosing $A= q^k=U^kS^k(V^k)^T $ and   $\delta A = \delta q^k $, we find $ \delta s^k / \delta q^k_{ij} =   (u_l^k  \otimes  v_l^k)_{ij}   $.  
The change of (\ref{J2}) by  perturbing    $q^k$ by $\delta q^k$ is thereby found       
\begin{eqnarray}
{\delta  J_2^k} = 2 %
 \sum_{l=1}^{r_k} \left(  -  s^{k}_l \frac { (u_l^k  \otimes  v_l^k)_{ij}    }{n^2_k}  + \frac{(s^{k}_l)^2}{(n^2_k)^2} q^k_{ij}  \right){\delta q^k_{ij} } 
  \label{gradSimple}  \,.
\end{eqnarray} 

The constraint (\ref{JRestrict}) can be incorporated by a simple consideration. 
A more formal derivation is found in appendix~\ref{formalGradient}.
The  variation  in one frame implies a variation  in the lab frame 
$ T^{\Delta^k} [{\delta J^k_2}/{\delta q^k}] $
which needs to be  redistributed to keep the constraint (\ref{JRestrict}) unaltered. 
Redistributing  to all  frames equally\footnote{
	The redistribution does not need to be equal for all frames to satisfy the constraint. 
		By using non-equal distributions, families of gradients can be constructed, a freedom which might be beneficial for the optimization, but which is not used in this report.}
 yields  
\begin{eqnarray}
\frac{\delta J_2}{\delta q^k}    =     \frac{\delta J^k_2}{\delta q^k}-  \frac 1 {\Ncm} T^{-\Delta^k} \left[  \sum_{k'=1}^{\Ncm}  T^{\Delta^{k'}} \left[ \frac{\delta J^{k'}_2}{\delta q^{k'}} \right] \right].
\label{distributeGrad}
\end{eqnarray} 
Any linear combination of gradients does not change the  constraint \ref{JRestrict}, so that any gradient method which constructs a solution 
from a start value satisfying the constraint and 
the linear space spanned by the gradients is guaranteed to satisfy  the constraint. 
Note, that the existence of an inverse is assumed in (\ref{distributeGrad}). 
Orthogonality is used in the formal proof of (\ref{distributeGrad})  in the appendix \ref{formalGradient}.     

Non-simple singular values were excluded in the derivation of  (\ref{gradSimple}).
Associated with such singular values are (left and right) linear spaces,  where any  orthogonal basis vector set serves as  singular vectors.
 If one assumes a double singular value,  $s_{r_k}^k =  s_{r_{k+1}}^k $,
only  the singular vectors associated with the  'first' singular value enters the gradient. 
This  non-uniqueness cannot be uniquely lifted in an infinitesimal neighborhood, since the different small perturbations of the $q^k$ field choose  different singular vectors \cite{Stewart1991}. 
We conclude  that if multiple singular values exist, $r_k$ should be chosen to either included or exclude all multiple singular values (for this optimization step) by which a unique gradient is defined. 
However,  using an ill-defined gradient during an intermediate step of an optimization might work without noticeable problems and  
no special measures where taken in the examples.          
 
For  $J_1$ the construction follows the same lines yielding  a   sub-gradient for each frame $J_1^k =\| q^k \|_* $,  
\begin{eqnarray}
{\delta  J_1^k}  =  \sum_{l=1}^{\bar d} \left(  u_l^k  \otimes  v_l^k   \right)_{ij} {\delta q^k_{ij} }\,  \label{gradSimpleJ1}, 
\end{eqnarray}   
where, $\bar d$, is the number of non-zero singular values. 
For vanishing singular values $J_1$ is non-smooth, so that only a   sub-gradient can be constructed. 
It is used with the same redistribution (\ref{distributeGrad}) to incorporate the constraint  without any further ado.

%
 \section{Uniqueness of the Decomposition and Extensions} 
 \label{secExtensions}
 
 We now proceed to  analyze two weaknesses  of  $J_2$ and $\bar J_2$. 
 If more modes per frame  or more frames are admitted than necessary  for a perfect decomposition, these extra degrees of freedom are not automatically zero. 
 By this, an essential and in practice disadvantageous non-uniqueness is introduced. 
 Both weaknesses 
 are shared with the previous approaches  \cite{ReissSchulzeSesterhennMehrmann2018,SchulzeReissMehrmann2019}. 
 We will find that, here, the  objective $J_1$ (\ref{J1}) has advantages. %
 \smallskip

 \paragraph{Many modes per frame}
 
 Consider first the limit $r_k\to d $ of $J_2$  (\ref{rankKJ}), where  $d=\min(m,n)$ is the maximal possible  rank of the  matrix $q \in \R^{m,n}$. %
 If all singular values are summed, i.e., a full rank description is used in the co-moving frames, \emph{any}  decomposition %
 yields $J_2  = 0 $ and, therefore,  no separation in co-moving structures is created.
 However, we have already argued  that the functional $J_1$ (\ref{J1}) acts at concentrating the norm to large singular values.
 It needs, in its basic form, the full spectrum of singular values. It  can therefore be used to define the limit $r_k\to d $.    %

 \paragraph{A  redundant frame}    
 A similar problem arises if more frames, a larger $N$, are used than needed for a low rank description.  
 In practical applications the necessary  number of frames might be unknown. 
 Ideally, the matrix $q^k$ of the frame $k$ which is not needed becomes identical to zero.    
 A simple example shows that this is in general not the case for $J_2$, by considering the acoustic example (\ref{qw2}) and adding a frame of zero velocity $\Delta(t) = 0$, i.e. the frame which 
 corresponds to  the standard  POD.
 A simple sine wave in space and time is of rank one in this POD frame, which could be interpreted as the lowest mode of a vibrating string. 
 If the time-frequency is chosen consistently with the transport velocity it  can also be decomposed into two traveling structures,
see \cite{ReissSchulzeSesterhennMehrmann2018} for a detailed example. 
 Thus, the standing wave can be compensated by extra terms in the moving frames so that the  extra frame induces  a non-trivial family of rank one solutions, including 
 a zero field in the POD frame. %
 Indeed, practical tests usually yield such  structures in  extra  frames, if  admitted. %
 
 We find that for our test cases the $J_1$ objective function removes such an  extra structure, yet slowly.
 It was shown  that can be rewritten as \emph{one} Schatten-1-norm, $J_1 = | Q|_* $, see (\ref{J1Schatten}). 
 The Schatten-1-norm was suggested as a heuristic method to reduce the rank \cite{FazelHindiBoyd2001}. 
 If this heuristic works it implies that a redundant frame should be set to zero, since by this    
 the rank of the block matrix $Q$ containing the various frames $q^k$ is reduced. 
 Therefore, we expect  that $J_1$ can remove redundant frames by driving the associated field to zero.

 \paragraph{Additional Regularizations} 
 
 The objective function $J_2$ (\ref{rankKJ}) can be improved by increasing the norm in each frame, by e.g. 
 adding a constant value everywhere and compensate it by the other frames. 
 This contradicts the intended purpose, but did rarely happen  in numerical tests.
 In such a case, a simple regularization by penalizing some norm in each frame can mitigate this issue. 
 
 The spatial or temporal modes are not necessarily  smooth. As an example consider two, equally shaped,  crossing pulses with different signs. 
 At the crossing moment they perfectly compensate, so that the time amplitude  in each frame is arbitrary and in fact  zero is typically produced by the optimization. 
 While mathematically sound, for physical models a continuous behavior in amplitudes is often assumed. 
 By adding a regularization penalizing, e.g. the second derivative, this smoothness of the amplitude could be promoted, however this  is not further investigated here.

\section{First Examples} 
\label{secFirstEx} 

Now we illustrate the gradient and the optimization  for some generic cases.
While the method allows for non-constant velocities and changing wave forms, we 
first consider the  case of two simple  waves (2w as annotation) with constant transport velocities $c_{1,2} = \pm 1 $ %
\begin{eqnarray}
q^{2\mathrm{w}}_{ij}=   r_1(x_i -c_{1} t_j  ) + r_2(x_i -c_{2} t_j  )\,.  
\label{qw2} 
\end{eqnarray}
The discrete field size is $m=100, n =50 $ in this section.  
The transported quantities 
$r_{i}(x)$ , $i=1,2$,  are Gaussian distributions around $L/4$ and $(3/4)L  $ with $ \sigma = 0.06\cdot L$, where the system length is $L=2\pi$.
While apparently simple, it yields a slow convergence with the  POD method \cite{ReissSchulzeSesterhennMehrmann2018} and 
 the two waves also hinder the use of a single transformation. 
With introducing shifts $\Delta^k=t\,c_k  $ and thereby  $T^{\Delta^k } [ r_k (x)]$ $ = r_k(x-c_k t   )$,  (\ref{qw2})  can be exactly represented in the form (\ref{decompGen}) with rank one   ($r_k=1$) and constant time amplitudes $\alpha^k$. 

\begin{figure}
	\begin{center}
	\begin{minipage}{.33\linewidth}
			\begin{center}
	        \includegraphics[width=\linewidth,clip=true,viewport=  0 0 390 270]{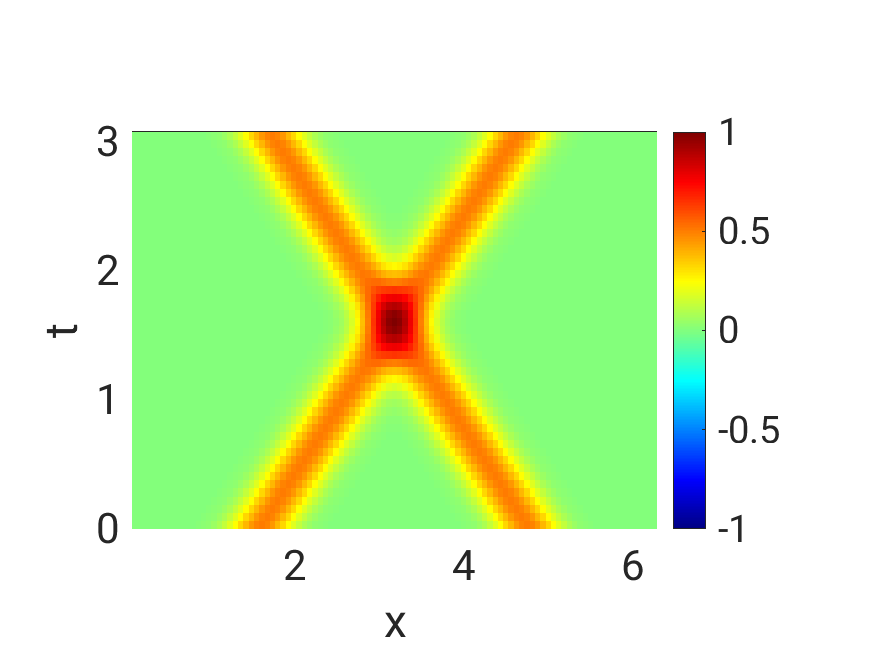}\\[1em]
			\includegraphics[width=\linewidth,viewport=  0 10 410 270]{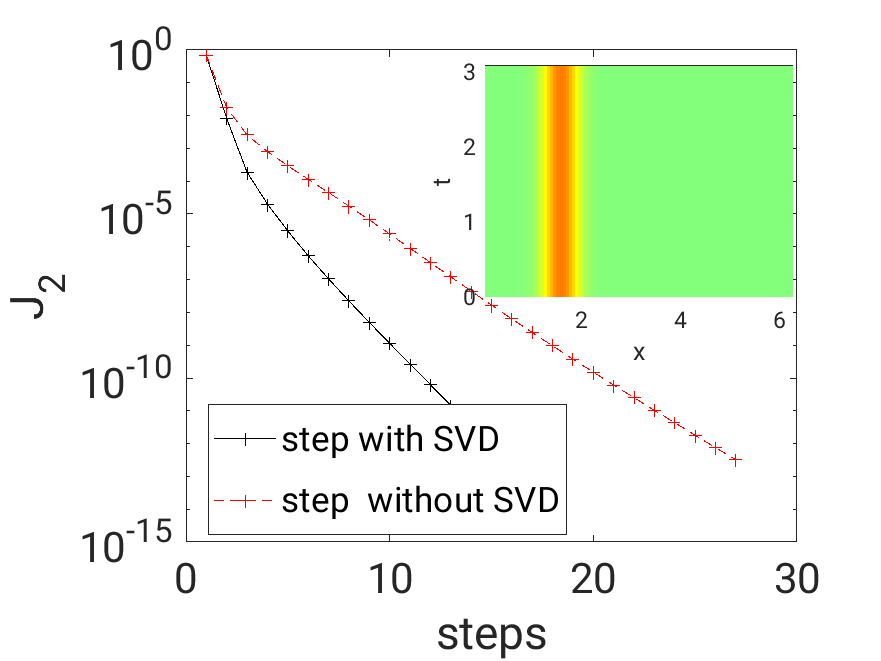}
			\end{center}
	\end{minipage}
	\begin{minipage}{.6\linewidth}	
			\includegraphics[width=\linewidth,clip=true,viewport=  20 10 390 300]{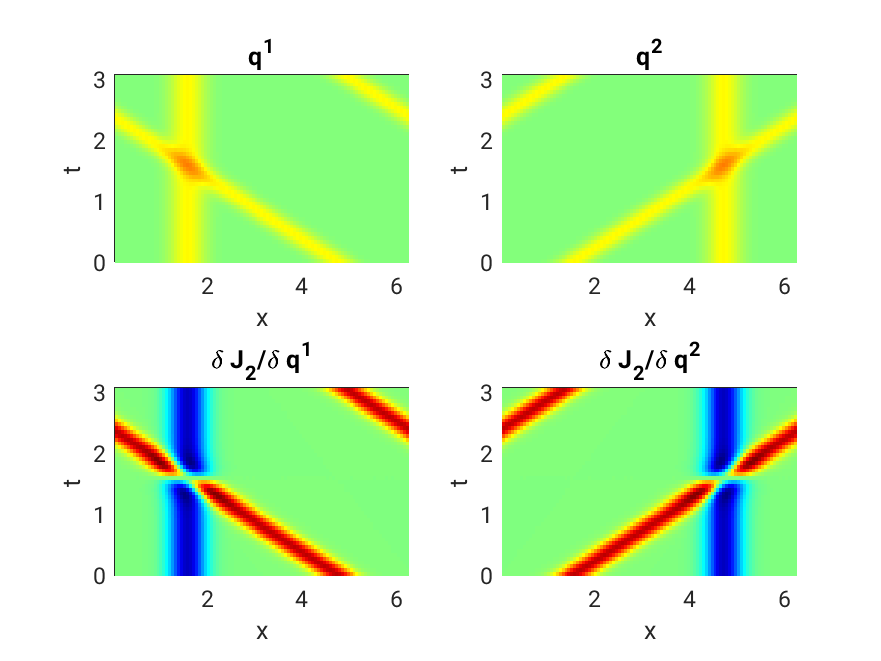}
	\end{minipage}
	\end{center}
	\caption{
		Top left: A generic case of two crossing waves in the space-time diagram, as the sum of  
		two simple waves (\ref{qw2}).
		Bottom left: The convergence  of  a  steepest descent method for different step estimators and the final result for the first frame, evidently of rank one.  
		Top middle and right:    
		The initial matrices in the two frames is the (shifted) equal distribution of the given data $q_k = T^{-\Delta^k} [q^{\mathrm{2w}}]/2 \quad k=1,2$. 
	    The final decomposition  of the test data $q^{\mathrm{2w}}$ based on  $ J_2$ with $r_k=1$  finds the expected simple standing pulse in each frame (left, middle).	    
	    Bottom middle and right: The gradient of the objective function $J_2$ (\ref{J2})  for the initial guess suggests to strengthen the vertical structure and remove the diagonal to 
	    improve the rank one approximation, see text.  
	    	}
	\label{fig_gardientRankOne}
	\label{fig_qg1}
	\label{fig_convergencepulse_oneMode}
\end{figure}

A method of convective decomposition should ideally find the structure  from the data $q^{\mathrm{2w}}$, which is  a vertical line in the co-moving frame,  evidently of rank one.
As an initial guess the original data $q$  is equally distributed over the two co-moving frames by transforming it to  each frame by $T^{-\Delta^k} $  divided by two,
see fig.\ (\ref{fig_gardientRankOne}), top left. 

The gradient  $(\delta J_2)/(\delta  q^k) $ is shown in fig.\ \ref{fig_gardientRankOne}, bottom middle and left. 
It suggests to increase  the wave at rest in each reference frame  (dark blue area) and to reduce the other wave (red). 
In the overlap region, little change is depicted. 
With an appropriate  step size as the diagonal structures are removed and re-attributed   to the respectively other frame.
By this one gradient step yields already a reasonable decomposition with only small errors in the overlap region. 

In spite of the high-dimensionality of the optimization, a simple steepest gradient method produces  a (numerically) perfect decomposition after a few iterations, see fig.~\ref{fig_convergencepulse_oneMode}, left bottom. 
The necessary  step size is estimated by  doubling test steps until  $J_2$ increases. The values of the  last three steps determine a parabolic interpolation, which yields the step as its minimum.  
For this estimator $J_2$ has to be evaluated several times, and with this the expensive SVD. 
The latter can be avoided by calculating the change of the singular values by 
\begin{eqnarray}
 \delta S = U^T \delta A V \label{overlap}
\end{eqnarray} 
where $\delta A$ is the change of the scaled gradient. 
Evaluating   $J_2$  on the base of this estimated singular values is much cheaper. 
It yields a slower convergence  (fig.~\ref{fig_convergencepulse_oneMode}, bottom left),
 but the overall cost to reach a certain quality where found to be cheaper for the cases studied in this report. 
 For this example, the step size estimation with exact singular values needs about three times as many  SVD evaluations.   
 \medskip

To test the objective function with higher ranks $r_k>1$ we consider the test function 
\begin{eqnarray}
	q^{\mathrm{2wd}}(x_i,t_j) &=&  p_1(x_i -c_{+} t_j  ) + t_j  \mu d_1(x_i -c_{+} t_j  )  \no \\
	                    &&    + p_2(x_i -c_{-} t_j  )   + t_j  \mu d_2(x_i -c_{-} t_j  )
	\label{qg2} 
\end{eqnarray}
where 
      $ p_\alpha(x) = \halb \exp\left( -\frac{(x-x_\alpha)^2}{\sigma^2} \right)   $ , 
      $ d_\alpha(x) =  \partial^2_x  p_\alpha(x ) $ %
      and $ c_\pm = \pm 1  $, $x_{1} = L/4$, $x_{1} = 3L/4$.  
The term proportional to $\mu$ mimics a diffusion with a diffusion  constant $\mu$. 
The data can be described by two modes in the two frames used before. 

\begin{figure}[h]
	\begin{center}
		\includegraphics[width=.4\linewidth,clip=true,viewport=0 0 400 310]{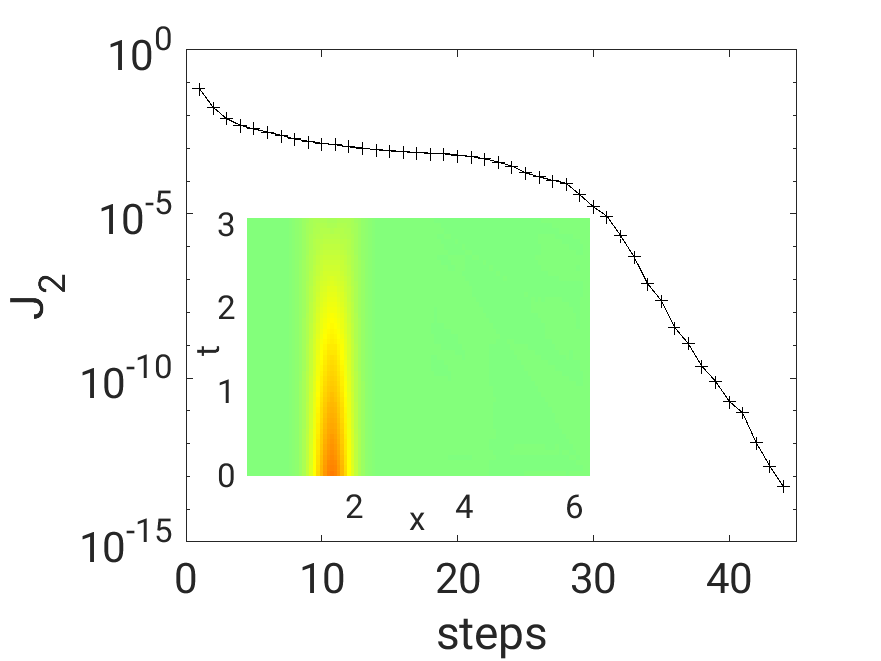}
	\end{center}
	\caption{%
		The final decomposition of the test data $q^{\mathrm{w2d}}$  with the objective functional $J_2$ using two modes in each frame ($r_k=2$).
		  The expected decaying structure of rank two is found  in each frame. 
		The convergence of the BFGS method is shown in the inset. }
	\label{fig_convergencepulse_twoMode}
\end{figure}

It turns out that the convergence rate of a steepest decent method is strongly reduced. 
This can be improved by starting with one mode in   (\ref{rankKJ}), $r_k=1$ and increasing $r_k$ when the value of  
  $ J_2$ saturates. 
  Instead, one can use a more sophisticated optimization method, which will be done in the following. 
  Here, we use the BFGS method implemented in the \textsc{Matlab} package \textsc{Hanso} described in 
   \cite{LewisOverton2013}. 
   The rank two description $r_k=2$  found by this is shown together with the convergence in fig.\ \ref{fig_convergencepulse_twoMode}.  
\smallskip

\begin{figure}
	\begin{center}
		\includegraphics[width=.4\linewidth]{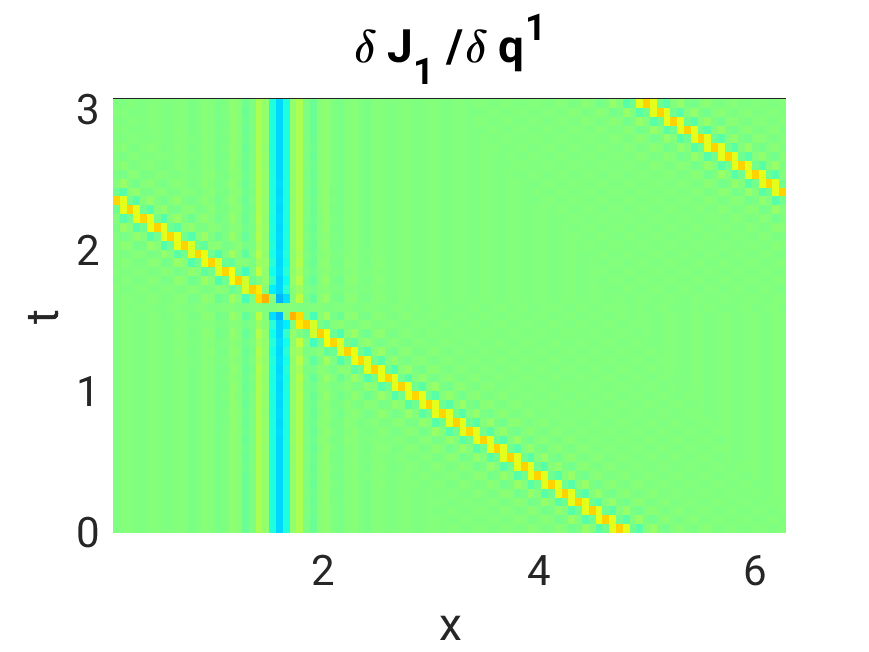}
		\includegraphics[width=.4\linewidth]{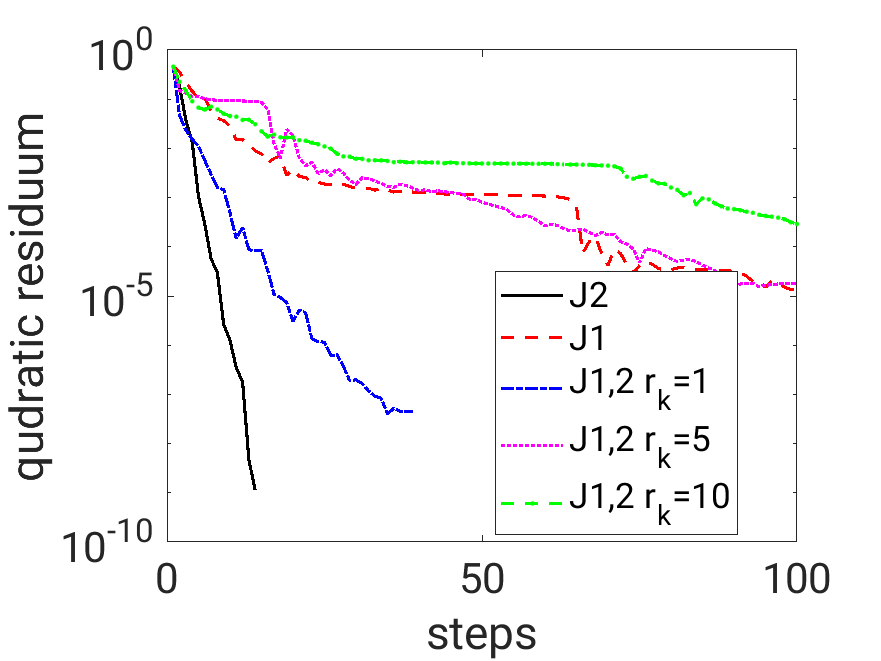}
	\end{center}
	\caption{Left: The gradient of the objective function $ J_1$ (\ref{J1})  for a equally split data $q_k =  T^{-\Delta^k} [q^{\mathrm{2w}}]/2 $.  
		     Right: The convergence of the $J_1$ and it approximations with a lower number of singular values (\ref{J1delta}) comparison to the $J_2$ functional.  }
	\label{fig_gardientJ1}
\end{figure}

Next, the performance of the $J_1$ objective function  (\ref{J1}) is investigated for the test case \ref{qw2}. 
 The gradient is shown in fig.~\ref{fig_gardientJ1}. 
Again the data $q^{g1}$ is equally distributed  over both frames.   
The gradient resembles  the one of the $J_2$ functional but the structures are smaller and a high frequency jitter is visible.  
This first gradient  cannot remove but only reduce the diagonal structures in the first step, which is reflected in a much worse convergence of the BFGS, see fig.~\ref{fig_gardientJ1}, right, where 
$    \sum (\tilde q_{ij} -q_{ij} )^2 / \sum(q_{ij})^2$ is plotted, with  $\tilde q$ constructed by (\ref{decompGen}) with one mode in each frame.

Beside the slower convergence, the $J_1$ function has, as noted before,  the disadvantage to demand all singular values. 
To overcome this, the modified objective function  (\ref{J1delta})  was derived. 
The convergence varies strongly  for different numbers  of  singular values, $r=1,5,10$. %
For one singular value it is close to the rate of $J_2$ which is, due to its similarity, not surprising.   
This suggests that adopting $r_k$ is a strategy to use $J_1$ with an improved convergence.

  \section{Boundary Treatment}

  \label{secBound}
  The shift operation implies that areas from outside the original or available domain are mapped inside. 
  These  unknown values have to be filled in some sensibly manner. 
  Extending the data by a constant induces, in general, jumps at the boundary which are harmful for a low rank description.   
  The derived formulation naturally extends to a boundary treatment, which  
  results in an implicit extrapolation for values outside the original domain. 
  \begin{figure}[!t]
  	
  	\includegraphics[width=.32\linewidth, viewport= 20 10 330 330, clip=true]{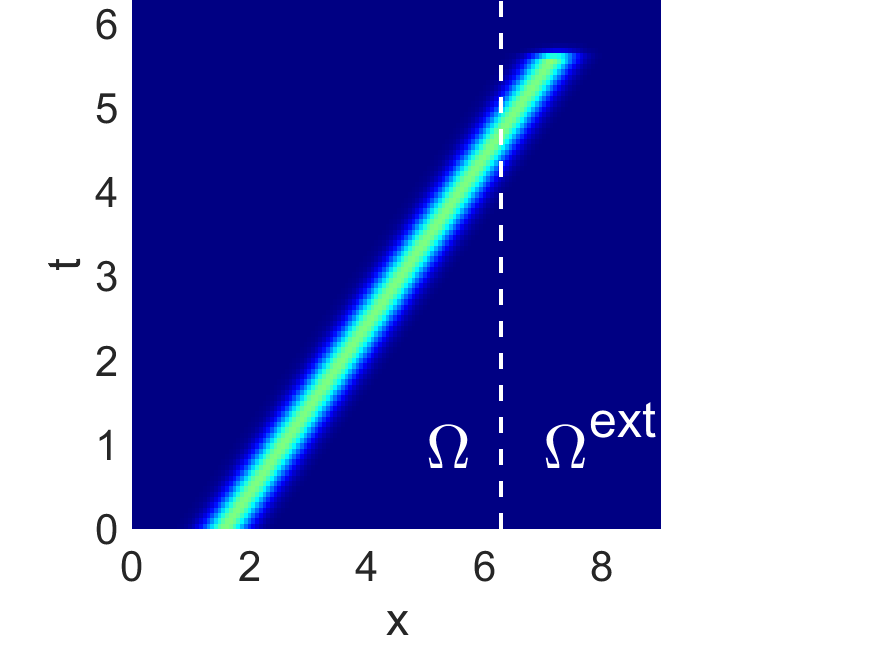}
  	\includegraphics[width=.32\linewidth, viewport= 20 10 330 330, clip=true]{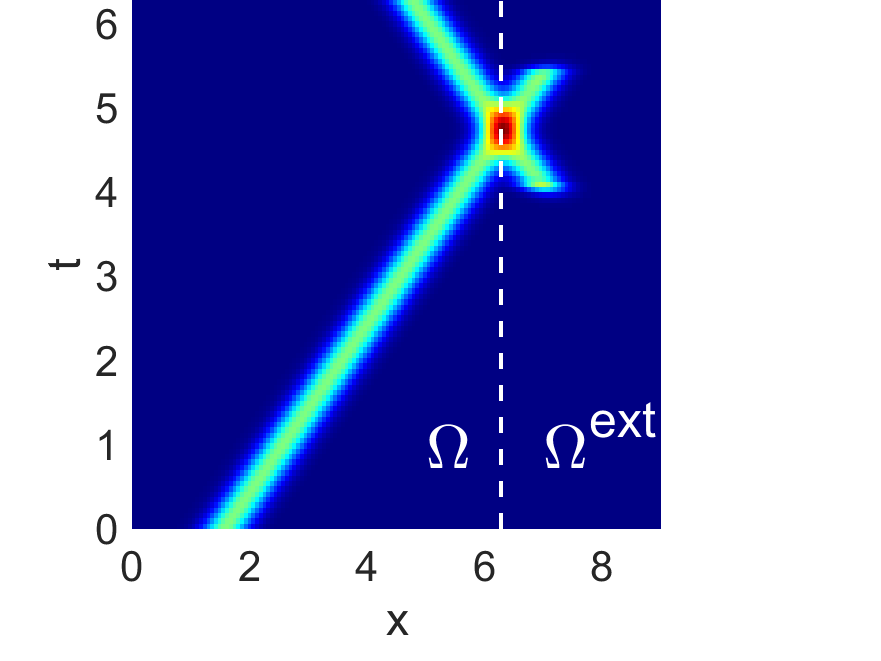}
  	\includegraphics[width=.32\linewidth, viewport= 20 10 330 330, clip=true]{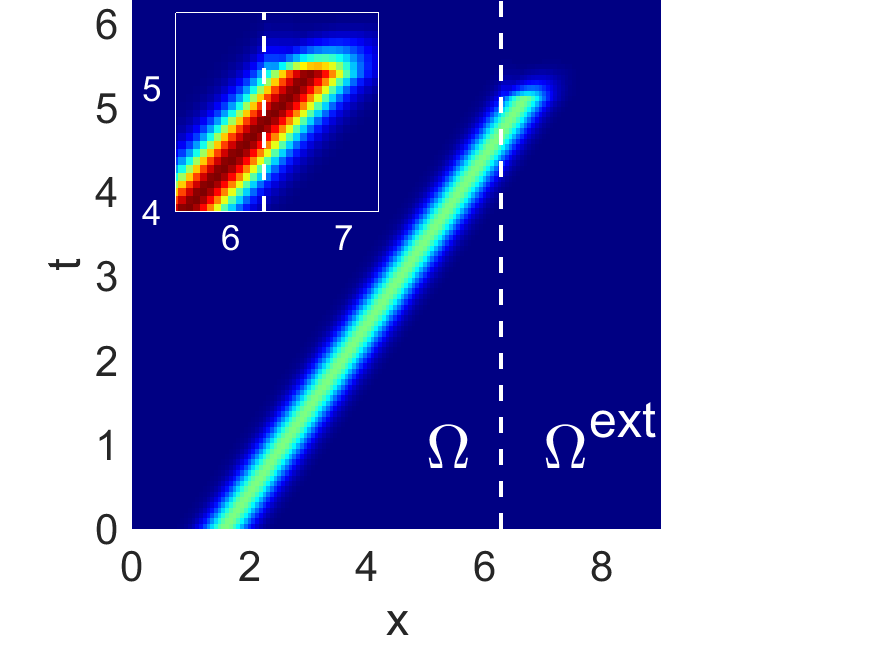}
  	\caption{
  		Left: A simple wave leaving the domain, domain boundary marked by the broken line. 
  		The values left of it are prescribed, the values right of it are constructed by the boundary treatment. 
  		Middle: A reflected wave. Here, the outgoing and the ingoing waves are implied by demanding  low rank in each co-moving frame.     
  		Right: The boundary extrapolation by the $J_1$ objective deviates slightly from rank one, visible in the detailed view in the inset, as an asymmetry of the extrapolated pulse at the boundary.
  		Color values scaled between 0 to 1, and 0 to 1/2 in the inset respectively. 
  	}
  	\label{fig_boundary}
  \end{figure}

  To handle boundaries we extend the original domain $\Omega$ to $\bar \Omega = \Omega \bigcup  \Omega^{\mathrm{ext}}$, in order to be big enough to hold all values that are mapped into the original domain by the shift operations, i.e. the spatial  domain is extended by the maximal   occurring shifts. 
  The values are initialized  by a constant, which is in principle arbitrary. 
  In this work,  we use typical values for all variables like the atmospheric pressure for  pressures or use a constant extrapolation, to give a good starting point for the optimization. 
  The boundary is treated simply by  relaxing the restriction (\ref{JRestrict}) to 
  \begin{eqnarray}
  w q &=& w \sum_{k=1}^{\Ncm }   T^{\Delta^k} [ q^k]  \label{JRestrictW},
  \end{eqnarray}
  where $w$ are weights defined by  
  \begin{eqnarray}
  w = \left\{
  \begin{array}{rl}
  	1 &  \mathrm{ in\;  } \Omega   \\
  	0 & \mathrm{ in\;  }  \Omega^{\mathrm{ext}} %
  \end{array}	
  \right.  
  \; .
  \end{eqnarray}
  By this, the decomposition restriction is enforced only in the original domain and
  the redistribution of the gradient (\ref{distributeGrad}) applies only  inside the original domain. %
  Values in $ \Omega^{\mathrm{ext}}$, outside the original domain, can be modified to improve a low rank approximation.

  To evaluate this procedure,  consider the simple example of  a wave leaving the domain.
  The wave is  Gaussian with $\sigma = 0.06 L $ with a domain length $L=2\pi$.   
  We use the objective function  $ J_2$ with only  one frame with the correct wave velocity of one. 
  The domain extension was initialized by zero.  
  The cut at the original domain boundary increases the rank in the co-moving frame, since this cut creates a diagonal structure in the co-moving frame. %
  It is obvious that extrapolating the solution along the transport  in the extended domain allows to  produce a rank one structure. 
  We find that this solution is indeed created, yet the extension has only  finite length. 
  As the Gaussian pulse  decreases exponentially  away from its center, the gain from extending the wave further into $ \Omega^{\mathrm{ext}}$  becomes very small, so that it is numerically a valid result. %

  The previous example  is extended by a second wave crossing at the boundary, which corresponds to a wave reflection, see  fig.~\ref{fig_boundary}, middle.
  Two frames with the corresponding velocities are used.  
  Again, the  extended  domain is filled with zeros. The desired extension of each wave is found. 
  It is noteworthy that this is structurally  the  boundary treatment by a mirror pulse  in the acoustic theory.

  If we allow the  two frames with velocities 1 and -1  for the simple, leaving wave,  we again face non-uniqueness, which can treated by a simple regularization if it poses a problem for the subsequent application.       
  Here for example,  adding  the 2-norm to the objective with a factor of $\epsilon = 10^{-4}$  on  all frames reduces the norm  of the  unneeded frame to 1/100 of the other frame after 200 BFGS steps. 
  
  The same strategy  works in principle also for the $J_1$ functional.  
  We find, however, a small, hardly visible kink at the boundary, see fig.~\ref{fig_boundary}, right,  which presumably originates from a competition between a quick singular value decay on the one hand and a small total norm on the other hand; 
  increasing the norm of the frame is in the $J_1$ in contrast to $J_2$ penalized and leads to more zeros in the modes. 
  However, we find  a factor of roughly 100 between the leading and the second singular value so that  in practice both objective functions, $J_1$    
  and $J_2$ yield a satisfactory  boundary treatment.

\section{Investigation of the Decomposition} 
\label{secInvestDecomp} 

In the following we analyze the gradient and the optimization based on this. 
The gradient of $J_2$ per frame, without incorporating the constraint (\ref{JRestrict}) is given by   (\ref{gradSimple}). 
By  factoring $ 1/n^2_k$ in (\ref{rankKJ})  the gradient can be rewritten as 
\begin{eqnarray}
\frac{\delta  J_2^k}{\delta q^k_{ij} } = %
&&
2 \left(  q^k_{ij} - \sum_l^{r_k}  s^{k}_l  (u_l^k  \otimes  v_l^k)_{ij} \right)  \frac 1 {n^2_k}     %
+
2 \left(   n^2_k  -  \sum_l^{r_k} (s^{k}_l)^2   \right) \frac{-q^k_{ij}} {(n^2_k)^2}  
\,.
\end{eqnarray} 
The first term is essentially the residuum of the truncated SVD 
\begin{eqnarray}
R^k =  q^k - \sum_l^{r_k}   s^{k}_l (u_l^k  \otimes  v_l^k)_{ij}  \,  \label{Rk}%
\end{eqnarray}
in each frame, which combines for  all frames  to the total residuum.
It simply states that it is advantageous to remove the residual part - the part not described by the rank $r_k$ approximation -  to improve the objective function.    
The modified gradient (\ref{distributeGrad}), which is constructed to respect the decomposition constraint (\ref{JRestrict}), redistributes all residua over all frames. 

The second term originates from the variation of the norm, so that the gradient produces the values  $q^k_{ij}$ as a factor in each frame.  
Decreasing the norm in one frame without changing the singular values reduces the respective term in the objective function,   
so that  %
it  can reduce the total objective function as long as the increase in other parts is not  overcompensating it. 
A change of the objective function by redistributing  the residuum is not helpful for finding a low rank  decomposition and is therefore  
collateral to the intention.         

This form  of the gradient yields a simple picture of the decomposition procedure.
The goal is to capture as much as possible by the  admitted number the singular vectors.  
The residuum  in each frame is the part which is likely  attributed to the wrong frame, as it might be possible included in  singular vectors in another frame.
Changing  along the gradient, which is constructed from  the residuum of all other frames, adds this parts to the other frames.
If there is an overlap of this gradient or residuum with the singular vectors of a frame, the singular values can be increased, improving the objective function.        
This allows to re-attribute parts  to other  frames. %
Thus, the singular values are increased if there is an overlap. For a finite step  the singular vectors are also modified.
The step estimator used in the section \ref{secFirstEx}  makes this interpretation explicit as the overlap of this residuum with the resolved (low rank) part in each frame is calculated.

The same interpretation holds for $\bar J_2$, eqn.~\ref{barJ2}. 
Here the gradient is simply ${\delta  J_2^k}/{\delta q^k_{ij} } = 2 R^k$ and a pure relocation  of the residuum between frames has less effect on the value of  the objective function.
This seems  favorable in principle, however, in numerical  examples little difference was observed between $J_2$  and $\bar J_2$.

\smallskip

The sub-gradient of the $J_1$ functional (\ref{gradSimpleJ1}) is simply the sum over all singular modes without the singular values. 
By this, the modes with small singular values are over-represented, so that using this as a gradient reduces these trailing modes strongest. 
This explains the way in which the separation is archived. 
But since the gradient  also contains the desired low rank part, it is less specific in separating  the desired and undesired parts, explaining its poorer performance for simple gradient based methods.

\section{Comparison with Former Formulations} 
\label{secFormer}  The connection with the two previous formulations, the shift\&reduce-method \cite{ReissSchulzeSesterhennMehrmann2018} and the 
residuum minimization method \cite{SchulzeReissMehrmann2019}, are now discussed.

The shift\&reduce-method 
builds on a shift of the data aligning with each co-moving frame, which  is followed  by a low rank approximation by truncating an  SVD  
\begin{eqnarray} 
P_{ r_k } T^{- \Delta^k} [q(x,t)]  ,
\end{eqnarray}
where $P_{ r }$ is a projector replacing a matrix by a low rank approximation by a trunctaed SVD with the leading $r$ singular values.  
It is acting as an effective filter for the structure, which is moving with the according frame, since  
 a transported structure changes, if at all, slowly in the co-moving frame
and  can be described well by the leading (few) modes. 
In contrast, structures moving relative to the frame are distributed over many modes and 
conversely have little weight in the first modes.  
Thus, the reduction in the co-moving frame is likely to be dominated  by the desired structure. 
The  interference by the other velocity component causes a pollution, for which an  iterative cleaning procedure was developed, where  
again, the  shifted\&reduced method is used on the residuum to identify possible corrections,    
 see \cite{ReissSchulzeSesterhennMehrmann2018} for details.

This can be compared with the $J_2$ or $\bar J_2$ optimization. 
Initially distributing $q$ equally over the frames is followed by calculating the leading singular values. 
The residuum $R^k$, (\ref{Rk}), is the essential part of the gradient.  	
Instead of doing a shift\&reduce of the residuum, which was observed to become less effective with each step as the co-moving structures become less and less pronounced in the residuum,  it here  acts as a gradient. 
Due to the constraint (\ref{JRestrict}), the residuum is redistributed over all frames (\ref{distributeGrad}). 
In a simple descent method, the gradient is added to the already identified low rank structure in order to increase  the leading singular values. 
The last point is evident from the step estimator which builds on the overlap of the gradient with the identified structure by 
setting $\delta A \sim \delta J_2/\delta q^k$  in  (\ref{overlap}). 
Calculating the overlap with the dominant structures helps to identify parts which are beneficial in this frame, which explains its quick convergence for simple cases.
Further, enforcing  the constraint in the gradient avoids  double accounting of structures, which are well represented in all frames, e.g. points in space time, which was a core problem in the report  \cite{ReissSchulzeSesterhennMehrmann2018}.

\smallskip 

The residuum was the target of the optimization in  the second method \cite{SchulzeReissMehrmann2019}. 
The objective function $\bar J_2$ (\ref{barJ2}) can  be used to bound  residuum\footnote{I am indebted to Philipp Schulze for pointing this out. } 
\begin{eqnarray}
R =  q-\tilde q 
= 
  \sum_{k=1}^{\Ncm} T^{\Delta^k}\!\left[q^k -\tilde q^{k}\right] \, , 
\end{eqnarray}
where  (\ref{JRestrict}) is used and the abreviation $\tilde q^k = \sum_{l=1}^{r_k} u^k_l s^k_l (v^k_l)^T $ for the low rank approximation in each frame.   
It follows for the Frobenius norm
\begin{eqnarray}
\|R\|_F& =& 
\left\|\sum_{k=1}^{\Ncm}   T^{\Delta^k}\!\left[q^k -\tilde q^{k} \right] \right\|_F 
\leq 
\sum_{k=1}^{\Ncm}   \left\| T^{\Delta^k}\!\left[q^k -\tilde q^{k} \right] \right\|_F \no\\
&=& 
 \sum_{k=1}^{\Ncm}  \left\| q^k -\tilde q^{k} \right\|_F 
 = 
\sum_{k=1}^{\Ncm}   \Bigl( n_k^2 -\sum_{l=1}^{r_k} (s^{k}_l)^2\Bigr)   = \bar J_2  \label{RleqJ2} 
\end{eqnarray} 
The linearity of  the shift operator was exploited and the conservation of the Frobenius norm, which is valid for the  (non-constant) shear transformations in space-time used here. 
Rotations in space-time would  also be admissible.    
By this, we  have basically the same optimization target as in the formulation of \cite{SchulzeReissMehrmann2019}, without directly enforcing a dyadic structure by using the entries of the dyadic pairs as optimization parameters as it is done there. 
In  contrast  $J_2$ cannot be used to bound the residuum, as the principally  arbitrary norm of each frame rescales the residuum of terms in (\ref{rankKJ}). 
If the norms are bound, as found in practice, a similar  estimate can be derived. 
\smallskip 

Since we do not have uniqueness of the solution for any of the methods we cannot guarantee that the final solutions are identical. 
However, if the data can exactly be decomposed into low rank parts $r_k \geq \mathrm{rank}(q^k)$, the objective functions
$J_2$ or $\bar J_2$ attain the minimal possible value zero.  
At the same time the residuum (\ref{residuum}) is zero. 
This is a fixed point of the shift\&reduce algorithm  and the goal of the residuum minimization so that such a decomposition is a solution for all these  methods.  
It is, however, not clear if the shift\&reduce algorithm is always able to reach this fixed point, 
as the convergence properties remained obscure in  \cite{ReissSchulzeSesterhennMehrmann2018}, despite its ability to find known decompositions in the considered cases.

The comparison with $J_1$ is more involved as it removes some non-uniqueness as discussed in the next section and, further, the rank minimization property is heuristic so that not always the minimal rank might be obtained, as seen in section \ref{secBound}.

\section{Application Example}
 \label{secPDC}
  
In this section, we apply the developed method to the numerical description of a pulse combustion chamber derived by data assimilation. 
The experiment and assimilation is presented  in  detail in \cite{GrayLemkeReissPaschereitSesterhennMoeck2017}. 
The situation is characterized by multiple transports with sharp fronts, including  transitional behavior hindering a composition into  steady state solutions  of the individual phenomena. 
This  data was  not tailored for the discussed method, but is a real world  example from engineering research.  
This application area is the original motivation of this method development, since the multiple sharp transports impede the use of classical mode based decompositions. 
The same configuration was investigated in \cite{SchulzeReissMehrmann2019} by the earlier mentioned residuum-based decomposition technique.  
\begin{figure}
	\includegraphics[width=.48\linewidth]{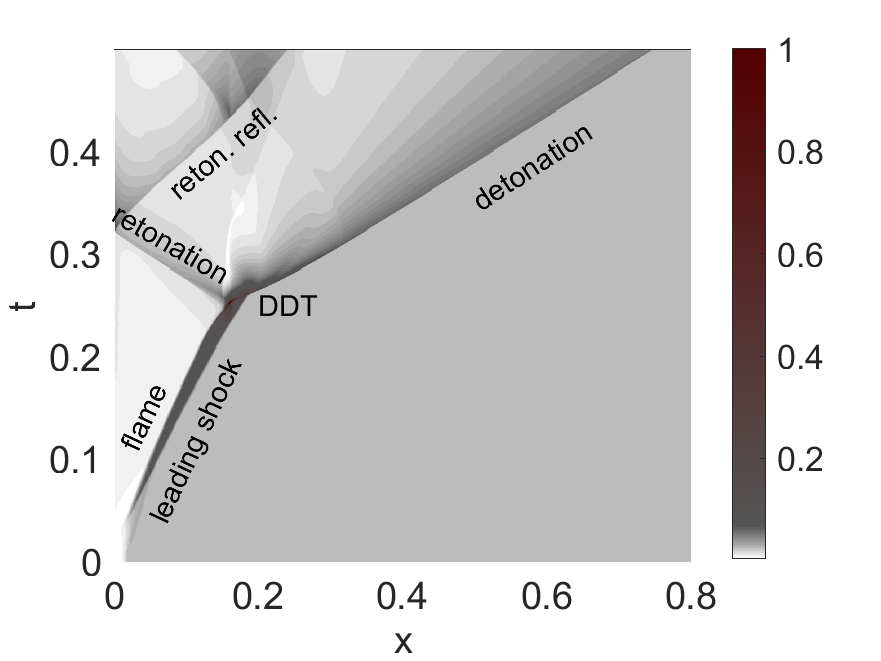}
	\includegraphics[width=.48\linewidth]{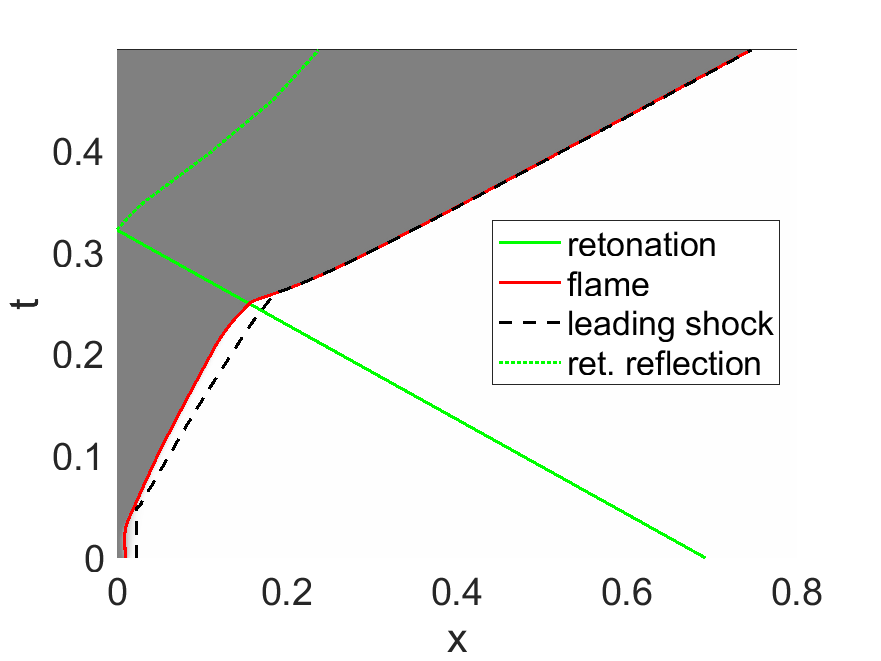}
	
	\caption{
		Left: The density field in space and time of the pulse detonation combustor shows a rich combination of transport phenomena. 
		     Close to $x=0$ a flame is ignited, which propagates and pushes like a piston creating a leading shock wave. 
		     Close to $x=0.2$ the flame changes to a detonation wave marked with DDT. 
		      This event is  accompanied by a high pressure peak creating a backwards traveling retonation  wave which is  reflected multiple times. 
		      The peak at the DDT event is larger than the other processes  by a factor of 10 .      
		      Right: The fraction of the reactant showing a burned ($Y=1$, gray) and unburned part($Y=0$, white). 
		      The flame front can be determined  by $Y=0.5$. 
		       The shock wave, found by a pressure threshold, merges the flame front after the DDT. The retonation wave location and its reflection are determined by inspection.    }
		    \label{fig_pdc_full}
\end{figure}

The original data is shown in figure \ref{fig_pdc_full}. 
It shows a traveling flame (deflagration)  which has, due to  the thermal expansion, a  piston-like effect and  thereby creates  a  shock wave with increasing strength. 
At a certain point, a  detonation wave is created from this flame, which is called a deflagration to detonation transition (DDT). 
First, this detonation has an increased detonation speed (overdriven detonation), which is then  continuously reduced to  a steady state detonation (CJ velocity). 
A high  pressure peak accompanies the DDT, also creating a backwards traveling shock called the retonation wave. %
This retonation wave is reflected at the left wall boundary. 
A varying cross-section of the combustion chamber, not depicted here, creates internal reflections at $x\approx 0.2$. 

The complexity challenges any data based decomposition approach. 
The data contains four variables: the density $\rho$, the momentum density $\rho u$ (density weighted velocity), the pressure $p$ and the fraction of the reactant $Y$, varying from $1$ (unburned) to $0$ (burned), see fig.~\ref{fig_pdc_full}, right. 
Each field is rescaled so that the maximum value is one, avoiding a units-dependent bias. 
All fields but $Y$ reached the maximum at the DDT location, where this peak excels all other structures by a factor of 10.
The resolution is 1024 points in space and 500 points in time.  
\begin{figure}
	\begin{minipage}{.48\linewidth}
		\begin{center}
			\small shifted POD approximation\\ 
			\includegraphics[width=\linewidth]{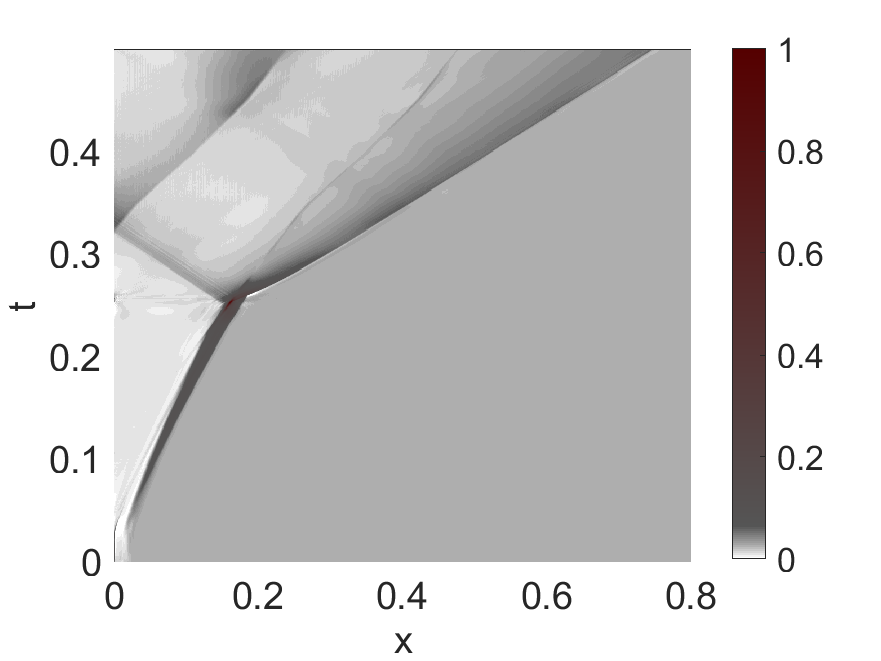}
		\end{center}
	\end{minipage}
	\begin{minipage}{.48\linewidth}
		\begin{center}
			\small error \\
			\includegraphics[width=\linewidth]{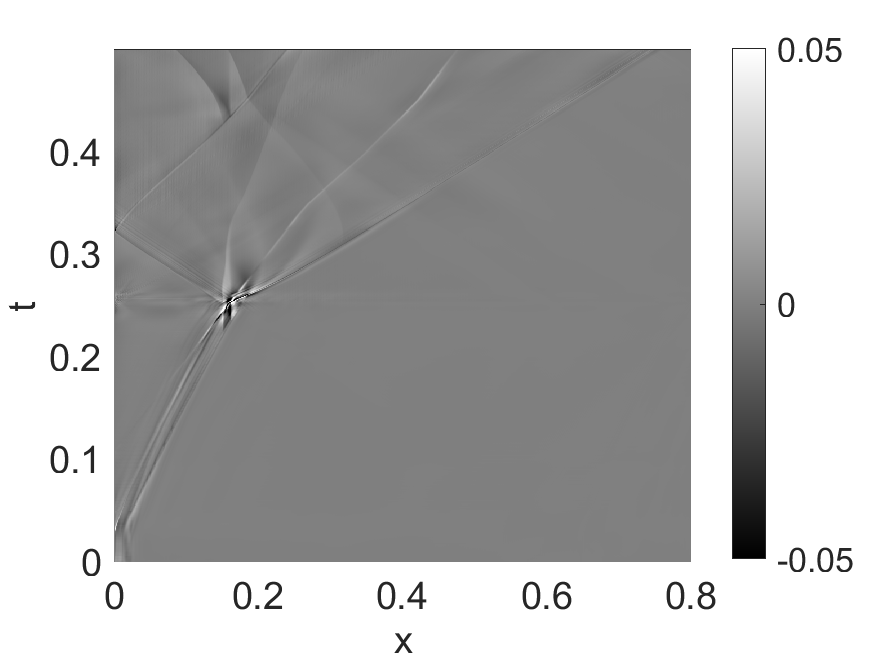}
		\end{center}
	\end{minipage}
	
	\begin{minipage}{.48\linewidth}
		\begin{center}
			\small POD\\
			\includegraphics[width=\linewidth]{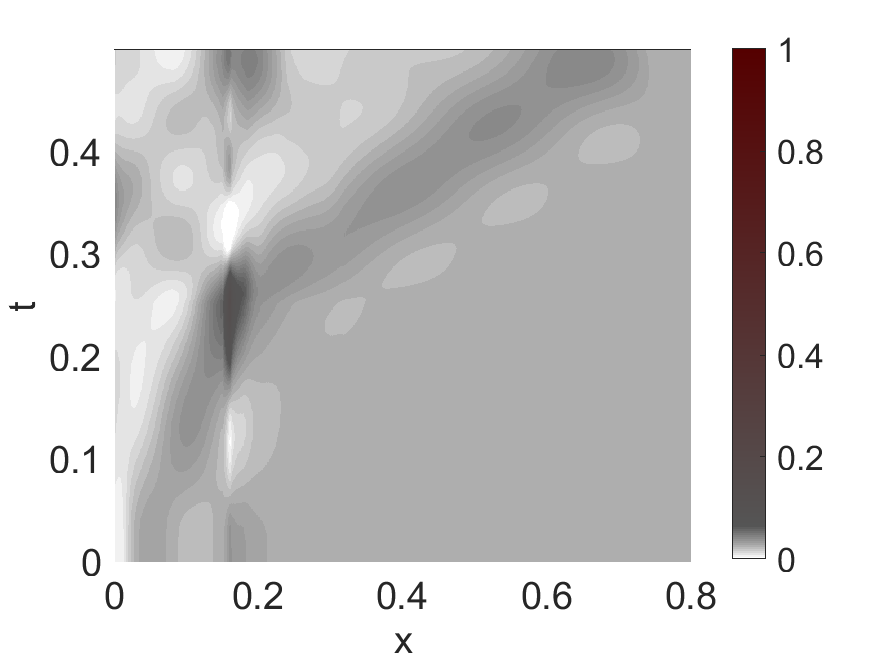}
		\end{center}
	\end{minipage}
	\begin{minipage}{.48\linewidth}
		\begin{center}
			\small error detail\\ 
			\includegraphics[width=\linewidth]{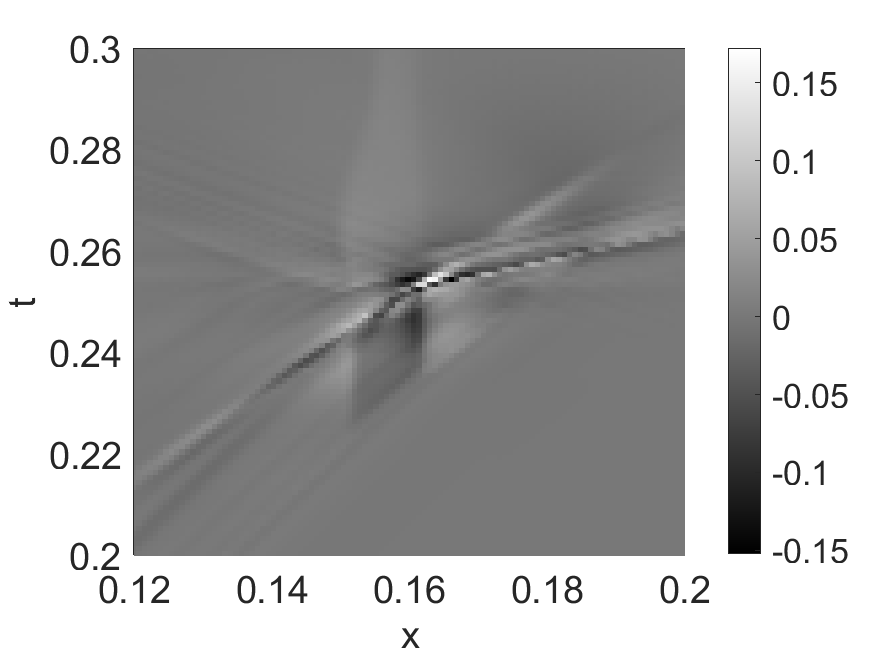}
		\end{center}
	\end{minipage}
	
	\caption{ Top left: The  approximation of the density by the shifted POD with seven modes shows a good agreement with the original data. 
		      Top right: The difference with the original data, (color-scale restricted) shows differences for the  DDT event, 
		              the flanks of the described waves and some missing  small, extra reflections. An additional wave is created parallel to the retonation reflection.    
		      Bottom left: Details of the difference close to the DDT event show the highest derivation.         
		      Bottom right: The POD with seven modes performs poorly.     }
	\label{fig_pdc_spod_pod}
\end{figure}

The velocity of the transported quantities is easily detected with the help of sharp structures. 
This  is done by peak search or threshold value search, whatever was found to give results matching visual inspection. 
Recently, an automatic velocity termination for the shifted POD was reported \cite{MendibleBruntonAravkinLowrieKutz2019}. 
Here, four velocity components are considered most important, shown in fig.~\ref{fig_pdc_full}, right:
(1) The flame or reaction front, as given by threshold of $Y=0.5$. 
(2) The shock wave created by this rapid expanding burning gas,   
combining after the DDT with the flame  to a detonation wave. 
The DDT event is associated with a very strong and narrow pressure peak creating a backward traveling shock wave called the retonation wave (3).
Its reflection (4) is also included in the decomposition approach. 
Secondary reflections of it close to $x=0.2$ and $t\approx0.45$ms are created by the carrying cross section of the combustion chamber (not depicted) and are not explicitly treated.         

The decomposition based on $J_2$ is done with four modes for the flame, and one for each of the other frames, in total 7 modes. 
To improve and accelerate the decomposition,  at first only the flame-frame  with one mode is permitted, so that the jump in $Y$ is directly attributed to this frame.  
This approximation is used to start a decomposition with three modes for the flame and one in all others, which again serves as a base with the finale  number of modes.
The BFGS method, in the  implementation \cite{LewisOverton2013}, is allowed to optimize for 50 steps using  25 modes to approximate the Hessian.

\begin{figure}
	\includegraphics[width=.48\linewidth]{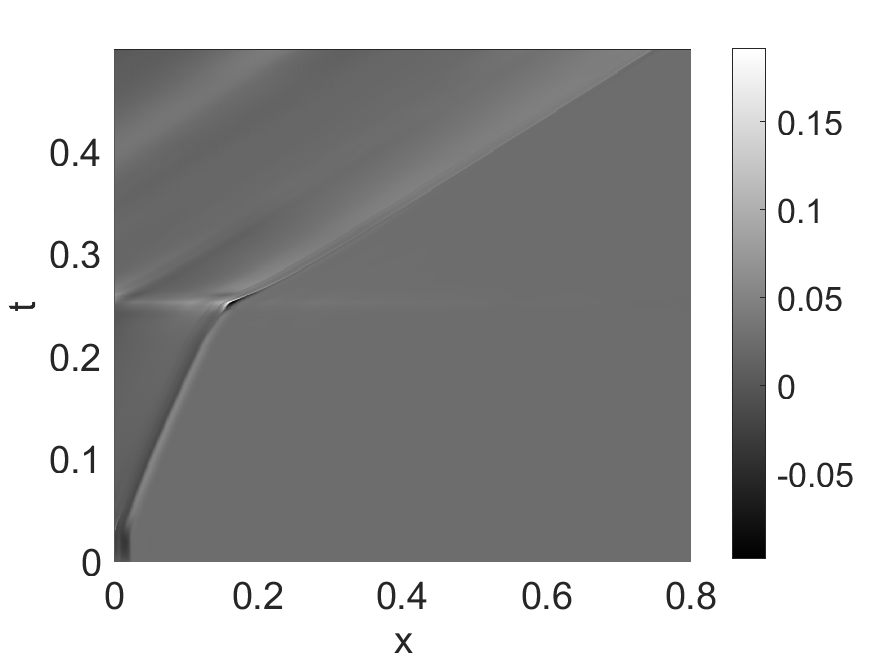}
	\includegraphics[width=.48\linewidth]{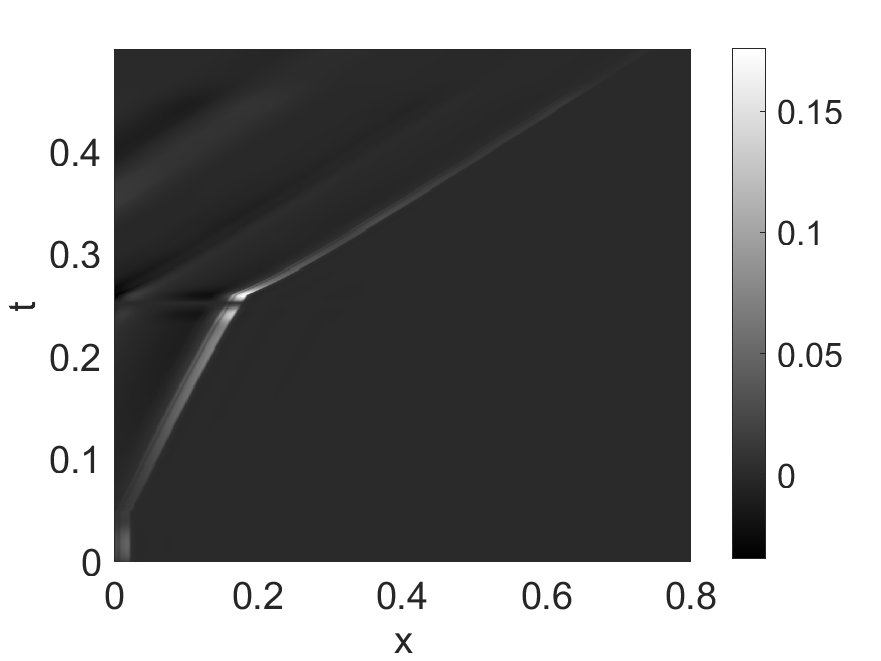}

	\includegraphics[width=.48\linewidth]{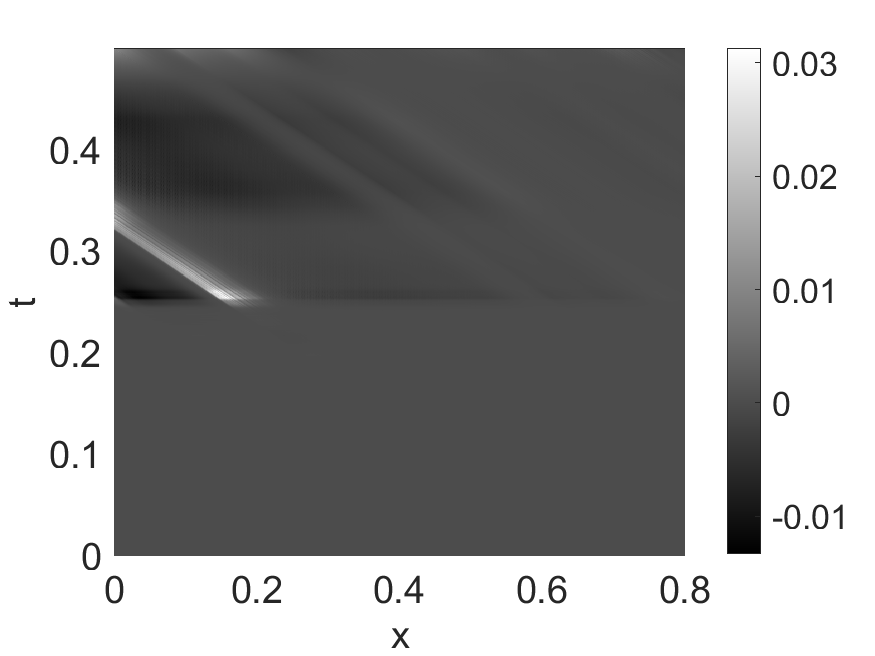}
	\includegraphics[width=.48\linewidth]{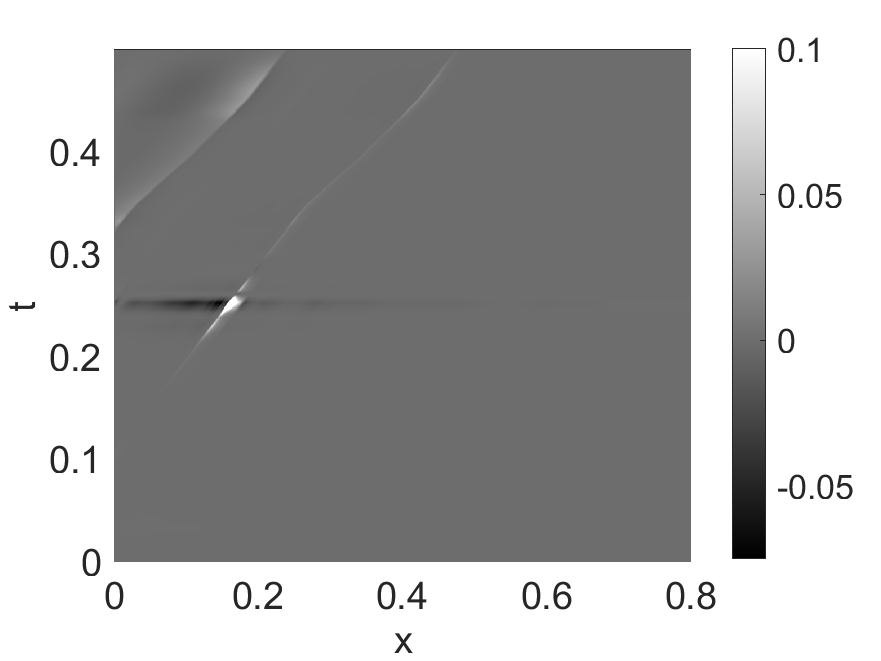}
	
	\caption{ The decomposition into the different frames, transformed to laboratory frame, density shown. While the flame (top left), the leading shock (top right) and the retonation wave (bottom left) 
		show the expected structure, the reflected retonation wave contributes strongly to the DDT peak (bottom right). The color scale has to be restricted to make the refelcted wave visible. 
        This is due to a approximative coincidence in the  velocity  of the reflected wave and  the  DDT event.      }
	\label{fig_pdc_frames}
\end{figure}
The final approximation is shown in fig.~\ref{fig_pdc_spod_pod}, top left. 
The main features  are well represented with a  relative error in Frobenius norm of $0.44\%$.
The difference to the original data is shown in    fig.~\ref{fig_pdc_spod_pod}, top right. 
Here, the color range is restricted to reveal fine structures.
The largest errors are close to the DDT location, presented  in the zoom in fig.~\ref{fig_pdc_spod_pod}, bottom right, which shows that the approximation fails to 
fully represent the maximal values there.   
It can be seen that the above mentioned reflections of the reflected retonation wave close to $x=0.2$ are missing. 
Further  an additional fine structure emerges from the DDT location parallel to the reflected retonation wave. 

While not perfect,  the qualitative and quantitative improvement is shown by a comparison with the POD, also with seven modes,  see fig.~\ref{fig_pdc_spod_pod} left bottom.  
The POD description has a   much higher $L_2$ error of $10.24\%$  and, what is considered  to be more important,   hardly shows any  features of the flow dynamics.

It is informative to investigate the individual co-moving frames, transformed into the lab frame for better inspection. 
The flame frame has its main features along the flame, as expected, see fig.~\ref{fig_pdc_frames}, top left.     
Some structures parallel to it in the top part, are visible, and a horizontal line of higher values at the time of the DDT event.  
The  shock frame, top right, looks similar, but with reduced values at times close to the DDT event. 
The retonation wave depicts mainly the desired retonation wave, bottom left. 
The reflected retonation wave does  not only show the desired wave, but it has an even much stronger contribution at the DDT location. 
By coincidence, the velocity of this wave is very similar to the DDT event so that from a pure decomposition perspective  it is profitable
to use the same frame to describe it. 
Here, due to a very similar velocity, the velocity criterion is not sufficient for a separation of physically different features.
In a practical  application one could add constrains, as zero weight before the reflection event, but we refrain here from doing so. 
The different horizontal lines visible close to the DDT event cancel, so that the original data is well represented.
It is difficult to determine whether this is optimal for the given data and modes, or if this is due to a not fully converged solution.  
However, this very high and narrow DDT peak was found to be challenging. 
It is so narrow that the data on this resolution is non-smooth, so that the  interpolation-based shift operation (here linear between the next points) has severe limitations in this region.  
The same holds true for the sharp gradients at the shocks. 
This technical difficulties currently hinder  a further improvement, but in light of the improvements over the POD, the  result is satisfactory  for this complex example.

\section{Conclusion}
\label{secConcl} 
A new way to derive  low rank approximation for transport phenomena was presented.
It is designed for configurations where multiple transports with different, possibly non-constant  transport speeds, are relevant. 
It builds on decomposing the original data into several frames,  where for each the transport is compensated by a  mapping compensating the transport. 
Each of these partial data can thereby be low rank, even if no low rank description for the original data is possible.
  
The method builds on the same decomposition concept as  previous works, the  the shift\&reduce shifted POD \cite{ReissSchulzeSesterhennMehrmann2018} and residuum based shifted POD \cite{SchulzeReissMehrmann2019}. 
It is, in contrast,  formulated in terms of  a singular values based objective functions, yielding a transparent mathematical interpretation. 
It contains a generic boundary treatment which produces an extrapolation of the values of the available domain outside the boundaries.   
Different possible objective functions are proposed and discussed.    
Namely, the functionals $J_2$ (\ref{rankKJ}) and $\bar J_2$ (\ref{barJ2}) built on comparing a truncated 2-norm of the singular values with the Frombnius norm. 
They differ in the weighting of the contribution per frame by the Frobenius norm associated with this frame.  
In contrast,  $J_1$ (\ref{J1}) builds on the sum of the singular values, the Schatten-1-norm. 
$J_1$ can be recast in a single Schatten-1-norm   (\ref{J1Schatten}), which connects the problem  to existing literature. 
The estimation of $J_1$ by a small number of singular values  (\ref{J1delta})  yields  a functional closely connected to $\bar J_2$, so that minimizing $\bar J_2$ for a small number of singular values in tendency improves $J_1$ alike. 
Moreover  $\bar J_2$ confines the residuum of a low rank description (\ref{RleqJ2})  and  thereby connects the  formulation  introduced here to the two existing  methods for the shifted POD.  

Numerical  tests underpin  an efficient numerical performance for generic and applied cases. 
Here, $J_2$ or $\bar J_2$ was found to converge in examples much faster  than $J_1$. 
This is however a preliminary finding as the special properties of $J_1$  as a Schatten-1-norm optimization were not exploited, i.e., no dedicated methods were utilized.

\section*{Acknowledgments}
I thank  Philipp Schulze, Philipp Krah  and Volker Mehr\-mann for valuable discussions.  
\section*{Funding}
{Funded by the Deutsche Forschungsgemeinschaft (DFG, German Research Foundation) - Projektnummer 200291049 - SFB 1029}
%

\appendix
\section{Formal calculation of the gradient}
\label{formalGradient}
To incorporate the constraint  (\ref{JRestrict}), an   auxiliary fields
$
\bar q^k(x,t)
$
without this constraint  can be defined. A projection 
\begin{eqnarray}
q^k  & = &  \bar q^k  + \frac 1  {\Ncm}   T^{-\Delta^k}  \left[ q - \sum_{k'=1}^{\Ncm }  T^{\Delta_{k'}}  \bar q^{k'}  \right]  
\end{eqnarray}
 is introduced to handle the constraint. 
 In the numerical implementation, the fields $\bar q^k$ are used to allow an  unconstrained optimization, 
and calculates $q^k $ at the end. Since  $q^k$  is a subset of $\bar q^k$ this does not reduce the set of  solutions.   
Linearizing this produces  
\begin{eqnarray}
\delta q^k  & = &  \delta \bar q^k  - \frac 1  {\Ncm}   T^{-\Delta^k}   \sum_{k'=1}^{\Ncm }  T^{\Delta_{k'}}  \delta \bar q^{k'} \no\\
&=&  \sum_{k'=1}^{\Ncm } \underbrace{\left( \delta_{k,k'}  - \frac 1  {\Ncm}   T^{-\Delta^k}     T^{\Delta_{k'}} \right)}_{=A^{k'k}}   \delta \bar q^{k'}\, .
\end{eqnarray}

It yields the  linear connection between $  \delta \bar q^{k'}_{\mathbf{i}'} = \sum_{k, \sI} A^{k',k}_{\sI',\sI } \delta q^k_\sI $, where the space-time index 
$\sI = (i,n)$ was introduced.   
This can simply be used with (\ref{gradSimple}) to provide 
\begin{eqnarray}
\frac{\delta  J_2}{\delta \bar q^k_\sI} 
= 
\sum_{k', \sI'}   \frac{\delta \tilde J_r}{\delta  q^{k'}_{\sI'}}  \frac{\partial q^{k'}_{\sI'} }{\partial \bar q^{k}_{\sI} }
= 
\sum_{k', \sI'}    (A^T)^{k,k'}_{\sI,\sI' } \frac{\delta  J_2}{\delta  q^{k'}_{\sI'}}
\end{eqnarray}     
or $({\delta  J_2}/\delta \bar q) = A^T  ({\delta  J_2}/{\delta \bar q})$.  
To calculate the transpose $A^T$  we assume that the transformation is an orthogonal matrix (and that $T^{-\Delta^k}$ is the inverse of $T^{\Delta^k} $ )  so that.
\begin{eqnarray}
(A^T)^{k',k}  &=& [\delta_{k,k'} - \frac 1  {\Ncm}        (T^{\Delta_{k}})^T (T^{-\Delta_{k'}})^T] \no \\
&=& [\delta_{k,k'} - \frac 1  {\Ncm}        (T^{-\Delta_{k}})^T (T^{\Delta_{k'}})^T]
\end{eqnarray} 
This is  true for the  time dependent  shift  considered here. 
Such a shift operator is a block diagonal matrix with simple shift matrices as blocks. %
Each of these simple shifts is an orthogonal matrix, which might be  approximate if interpolation between grid points is used. 
But also rotations in space time, as used in  \cite{SesterhennShahirpour2016},  fulfill this property.

\bibliographystyle{plain}
\bibliography{local}
\end{document}